\documentclass[12pt]{article}
\usepackage[utf8]{inputenc}
\usepackage[T1]{fontenc}
\usepackage{amssymb}
\usepackage{amsfonts}
\usepackage{floatrow}
\usepackage{amsthm}
\usepackage{graphicx}
\usepackage{amsmath}
\usepackage{eucal}
\usepackage{stmaryrd}
\usepackage{caption}
\usepackage{textcomp}
\usepackage{mathtools}
\usepackage{amsfonts,epsfig,epstopdf,titling,url,array}
\numberwithin{equation}{section}

\newtheorem{prop}{Proposition}

\newcommand{\beq}{\begin{equation}}
\newcommand{\RN}[1]{%
  \textup{\uppercase\expandafter{\romannumeral#1}}%
}
\newcommand{\eeq}{\end{equation}}

\newcommand{\beqs}{\begin{equation*}}
\newcommand{\eeqs}{\end{equation*}}

\usepackage[
top    = 2.75cm,
bottom = 2.50cm,
left   = 2.3cm,
right  = 2.3cm]{geometry}

\title{About a Proposition on Escobar\textquotesingle s paper ``The Geometry of the First Non-zero Stekloff Eigenvalue''}  

\author{
\small Leoncio Rodriguez Quiñones, email: rodriguez-l@javeriana.edu.co\\ \small Departamento de Matem\'{a}ticas\\
\small Pontificia Universidad Javeriana, Bogot\'{a}, Colombia
}

\begin{document}
\bibliographystyle{plain}
\large{

\date{}
\maketitle

\begin{abstract}
Let $(M^{2},g_{0})$ be a compact manifold with boundary, and let $g$ and $g_{0}$ be conformally related by $g=e^{2f}g_{0}$. We show that the inequality
$$\nu_{1}(g)\geq\Big(\max_{x\in\partial M}e^{-f(x)}\Big)\nu_{1}(g_{0})$$ stated in Proposition 2 in \cite{Escobar1}, is only possible when the equality is achieved. In order to achieve such equality, it is required that the function $f$ be constant on $\partial M$, as it is mentioned in \textit{Remark 3} also in \cite{Escobar1}. Hence, the scope of this inequality is less broad than the one suggested by the Proposition.

\end{abstract}

Let $(M^{2},g)$ be a compact Riemannian manifold with boundary of dimension $2$. In this brief note, we make a comment concerning a result obtained by Escobar in \cite{Escobar1}. We consider the boundary value problem given by
\begin{equation}\label{stekloff}
\begin{aligned}
    \Delta u &= 0 \ \  \text{\ \  \ in } M,\\
    \frac{\partial u}{\partial n_{g}} &= \nu u \ \  \text{ on } \partial M.
\end{aligned}    
\end{equation}
Where $\partial u/\partial n_{g}$ denotes the derivative of $u$ in the direction of the outward unit normal vector on $\partial M$. The boundary value problem (\ref{stekloff}) is the so-called \textit{Stekloff eigenvalue problem}, introduced by \cite{Steklov} in $1902$, and the $\nu$\textquotesingle s are known as the \textit{Stekloff eigenvalues}.  

Here we are interested in making a comment about a lower bound for the first non-zero Stekloff eigenvalue denoted by $\nu_{1}$. Namely, we look at the following inequality
\begin{equation}\label{stekloff1}
\nu_{1}(g)\geq\Big(\max_{x\in\partial M}e^{-f(x)}\Big)\nu_{1}(g_{0}).
\end{equation}
Inequality (\ref{stekloff1}) is stated in Proposition $2$ in \cite{Escobar1}, and the metrics $g$ and $g_{0}$ are conformally related, i.e, $g=e^{2f}g_{0}$, where $f\in C^{\infty}(M)$. This inequality is used by Escobar to prove Proposition $4$ in \cite{Escobar1}; it is also used in three examples presented by the author in order to show how inequality  (\ref{stekloff1}) works. 

In what follows, we show that the scope of inequality (\ref{stekloff1}) is not as broad as the the result suggests. Namely, we show that this inequality cannot be strict and that in fact, it is an equality which only occurs if the function $f$ is constant on $\partial M$ as \textit{Remark} $3$ in \cite{Escobar1} points out.

Let us being with $\phi\in C^{1}(\overline{M})$ and consider the \textit{Rayleigh Quotient} of $\phi$ denoted by $Q_{g}$ and defined by (see, \cite{Escobar1})
\begin{equation}\label{stekloff2}
  Q_{g}(\phi) = \frac{\int_{M}|\nabla \phi|_{g}^{2}dv_{g}}{\int_{\partial M}\phi^{2}d\sigma_{g}}.
\end{equation}
We argue that the inequality $Q_{g}(\phi) \geq\Big(\displaystyle\max_{\partial M}e^{-f}\Big)Q_{g_{0}}(\phi)$ cannot be strict; this implies that inequality (\ref{stekloff1}) cannot be strict either, and that in order for it to be true, it is required that the function $f$ be constant on $\partial M$.

\begin{prop}
    Inequality $Q_{g}(\phi) \geq\Big(\displaystyle\max_{\partial M}e^{-f}\Big)Q_{g_{0}}(\phi)$ only holds for the equality case, and such equality occurs if $f$ is constant on $\partial M$. If $f$ is not constant, equality is not guaranteed.
\end{prop} 

\begin{proof}
    First, we have that $g$, $g_{0}$ are conformally related. That is $g=e^{2f}g_{0}$, where $f\in C^{\infty}(M)$. Notice that in local coordinates, $dv_{g}=\sqrt{|\det(g_{ij})|}dx^{1}\wedge\dots\wedge dx^{n}$. Since $\det(g_{ij})=e^{2nf}\det(g_{0ij})$, we have $dv_{g}=e^{nf}\sqrt{|\det(g_{0ij})|}dx^{1}\wedge\dots\wedge dx^{n}=e^{nf}dv_{g_{0}}$.
    On the other hand, $|\nabla_{g}\phi|^{2}_{g}=g(\nabla_{g}\phi,\nabla_{g}\phi)$ and $\nabla_{g}\phi=g^{ij}\partial_{i}\phi\partial_{j}=e^{-2f}g^{ij}_{0}\partial_{i}\phi\partial_{j}=e^{-2f}\nabla_{g_{0}}\phi$ (here we are using the notation $\partial_{i}=\frac{\partial}{\partial x^{i}}$). Then, $|\nabla_{g}\phi|^{2}_{g}=e^{2f}g_{0}(e^{-2f}\nabla_{g_{0}}\phi,e^{-2f}\nabla_{g_{0}}\phi)=e^{-2f}g_{0}(\nabla_{g_{0}}\phi,\nabla_{g_{0}}\phi)=e^{-2f}|\nabla_{g_{0}}\phi|^{2}_{g_{0}}$. Thus, 
    \begin{equation}
        \int_{M}|\nabla\phi|^{2}_{g}dv_{g}=\int_{M}|\nabla\phi|^{2}_{g_{0}}e^{(n-2)f}dv_{g_{0}}.
    \end{equation}
    Notice also that $d\sigma_{g}=e^{(n-1)f}d\sigma_{g_{0}}.$ Then, if $n=2$ we get
    \begin{equation}
        Q_{g}(\phi)=\frac{\int_{M}|\nabla\phi|^{2}_{g}dv_{g}}{\int_{\partial M}\phi^{2}d\sigma_{g}}=\frac{\int_{M}|\nabla\phi|^{2}_{g_{0}}dv_{g_{0}}}{\int_{\partial M}\phi^{2}e^{f}d\sigma_{g_{0}}}.
    \end{equation}
    
    Suppose 
    \begin{equation}\label{stekloff3}
        \frac{\int_{M}|\nabla\phi|^{2}_{g_{0}}dv_{g_{0}}}{\int_{\partial M}\phi^{2}e^{f}d\sigma_{g_{0}}}\geq\Big(\displaystyle\max_{\partial M}e^{-f}\Big)\frac{\int_{M}|\nabla\phi|^{2}_{g_{0}}dv_{g_{0}}}{\int_{\partial M}\phi^{2}d\sigma_{g_{0}}}
    \end{equation}
    holds. Then, since we are not dealing with the zero eigenvalue, we have that $\nabla\phi$ is not identically zero. Then, (\ref{stekloff3}) implies
    \begin{equation}
        \int_{\partial M}\phi^{2}(1-Ke^{f})d\sigma_{g_{0}}\geq 0, \text{ with } K=\displaystyle\max_{x\in\partial M}e^{-f(x)}.
    \end{equation}
Set $A^{+}=\{x\in\partial M|\phi^{2}(1-Ke^{f(x)})>0\}$, $A^{-}=\{x\in\partial M|\phi^{2}(1-Ke^{f(x)})\leq 0\}$ and for a given function $h:\overline{M}\rightarrow \mathbb{R}$ denote its positive and negative part by $h^{+}(x)=\max\{h(x),0\}$ and $h^{-}(x)=\max\{-h(x),0\}$ respectively. Then,
\begin{equation}\label{stekloff4}
\begin{aligned}
0\leq\int_{\partial M}\phi^{2}(1-Ke^{f})d\sigma_{g_{0}} &=\int_{\partial M}\phi^{2}(1-Ke^{f})^{+}d\sigma_{g_{0}}-\int_{\partial M}\phi^{2}(1-Ke^{f})^{-}d\sigma_{g_{0}}\\
&=\int_{A^{+}}\phi^{2}(1-Ke^{f})d\sigma_{g_{0}}+\int_{A^{-}}\phi^{2}(1-Ke^{f})d\sigma_{g_{0}}.
\end{aligned}
\end{equation}
Therefore, from (\ref{stekloff4}) we conclude
\begin{equation}\label{stekloff5}
    \int_{A^{-}}\phi^{2}(Ke^{f}-1)d\sigma_{g_{0}}\leq\int_{A^{+}}\phi^{2}(1-Ke^{f})d\sigma_{g_{0}}.
\end{equation}
Notice $1=e^{f(x)}e^{-f(x)}\leq Ke^{f(x)}$. Thus, $1-Ke^{f(x)}\leq 0$ for all $x\in\partial M$, which implies $A^{+}$ is the empty set, and $Ke^{f(x)}-1\geq 0$ for all $x\in\partial M$ in particular $Ke^{f(x)}-1\geq 0$ on $A^{-}$. Then, $\phi^{2}(Ke^{f(x)}-1)\geq 0$ on $A^{-}$, and by monotonicity of the integral we have $\int_{A^{-}}\phi^{2}(Ke^{f}-1)d\sigma_{g_{0}}\geq 0$. Using (\ref{stekloff5}) we obtain
\begin{equation}\label{stekloff6}
    0\leq\int_{A^{-}}\phi^{2}(Ke^{f}-1)d\sigma_{g_{0}}\leq\int_{A^{+}}\phi^{2}(1-Ke^{f})d\sigma_{g_{0}}=0.
\end{equation}
Therefore, 
\begin{equation}\label{stekloff7}
    \int_{A^{-}}\phi^{2}(Ke^{f}-1)d\sigma_{g_{0}}=0.
\end{equation}
This implies that
\begin{equation}\label{stekloff8}
    \int_{\partial M}\phi^{2}(Ke^{f}-1)d\sigma_{g_{0}}=0.
\end{equation}
We claim that $Ke^{f(x)}=1$ on $\partial M$. 

Let $\mathcal{A}=\{(\mathcal{U},x)|x:\mathcal{U}\rightarrow \mathbb{R}\}$ be an oriented atlas of $\partial M$. By indexing $\mathcal{U}=\mathcal{U_{\alpha}}$ for some set of indices $A$, we can write $\{U_{\alpha}\}_{\alpha\in A}$; this is an open cover of $\partial M$ made up of chart domains. By Theorem $2.23$ in \cite{Lee}, there exists a partition of unity $\{\psi_{\alpha}\}_{\alpha\in A}$ subordinate to $\{\mathcal{U}_{\alpha}\}_{\alpha\in A}$. Since $\partial M$ is compact, there is $\{\mathcal{U}_{\alpha}\}^{N}_{\alpha=1}$, such that $\partial M\subseteq\bigcup^{N}_{\alpha=1}\mathcal{U}_{\alpha}$. Then, by definition of partition of unity, $supp\text{ }\psi_{\alpha}\subseteq\mathcal{U}_{\alpha}$, and we get
\begin{equation}\label{stekloff9}
    \int_{\partial M}\phi^{2}(Ke^{f}-1)d\sigma_{g_{0}}=\sum_{\alpha\in A}\int_{\partial M}\psi_{\alpha}\phi^{2}(Ke^{f}-1)d\sigma_{g_{0}}=\sum_{\alpha\in A}\int_{\mathcal{U}_{\alpha}}\psi_{\alpha}\phi^{2}(Ke^{f}-1)d\sigma_{g_{0}}.
\end{equation}

Using integration over a Riemannian manifold we have
\begin{equation}\label{stekloff10}
    \int_{\mathcal{U_{\alpha}}}\psi_{\alpha}\phi^{2}(Ke^{f}-1)d\sigma_{g_{0}}=\int_{x(\mathcal{U}_{\alpha})}\psi_{\alpha}(x)\phi^{2}(x)(Ke^{f(x)}-1)\sqrt{|\det(g_{ij})|}dx.
\end{equation}
Now, since $1\leq Ke^{f(x)}$ for all $x\in\partial M$, we have that $\phi^{2}(x)(Ke^{f(x)}-1)\geq 0$. Also, notice that since $\psi_{\alpha}\geq 0$ then,  $\psi_{\alpha}(x)\phi^{2}(x)(Ke^{f(x)}-1)\geq 0$ on $\mathcal{U}_{\alpha}$. This implies 

\begin{equation}\label{stekloff11}
   \int_{x(\mathcal{U}_{\alpha})}\psi_{\alpha}(x)\phi^{2}(x)(Ke^{f(x)}-1)\sqrt{|\det(g_{ij})|}dx \geq 0.
\end{equation}

We have showed that $\int_{\partial M}\phi^{2}(Ke^{f}-1)d\sigma_{g_{0}}=0$. Then, 
\begin{equation}\label{stekloff12}
    \begin{aligned}
        0 &= \int_{\partial M}\phi^{2}(Ke^{f}-1)d\sigma_{g_{0}}\\
        &= \sum_{\alpha\in A}\int_{\partial M}\psi_{\alpha}\phi^{2}(Ke^{f}-1)d\sigma_{g_{0}}\\
        &= \sum_{\alpha\in A}\int_{\mathcal{U}_{\alpha}}\psi_{\alpha}\phi^{2}(Ke^{f}-1)d\sigma_{g_{0}}\\
        &= \sum_{\alpha\in A}\int_{x(\mathcal{U}_{\alpha})}\psi_{\alpha}(x)\phi^{2}(x)(Ke^{f(x)}-1)\sqrt{|\det(g_{ij})|}dx.
    \end{aligned}
\end{equation}
Since (\ref{stekloff11}) holds for all $\alpha\in A$. We have that
\begin{equation}\label{stekloff13}
    \int_{x(\mathcal{U}_{\alpha})}\psi_{\alpha}(x)\phi^{2}(x)(Ke^{f(x)}-1)\sqrt{|\det(g_{ij})|}dx = 0
\end{equation}
for all $\alpha\in A$ as well. Since $\phi$ is arbitrary, $\psi_{\alpha}\phi$ is given; also  $supp\text{ }\psi_{\alpha}\phi$ is contained in $ \mathcal{U}_{\alpha}$ and it is compact. Then, by the Fundamental Lemma of Calculus of Variations (see, page 22 in \cite{gelfand}), we must have $Ke^{f(x)}-1=0$ almost everywhere on $\mathcal{U}_{\alpha}$; by continuity of $f$, this implies that $Ke^{f(x)}=1$ everywhere on $\mathcal{U}_{\alpha}$ for any $\alpha=1,2,\dots,N$. Since $\partial M\subseteq \bigcup_{\alpha=1}^{N}\mathcal{U}_{\alpha}$, then $Ke^{f(x)}=1$ on $\partial M$, and $f(x)=\ln{(1/K)}$ for all $x\in\partial M$, i.e, $f$ is constant.
Notice also then, that (\ref{stekloff8}) rules out  $\int_{\partial_{M}}\phi^{2}(1-Ke^{f})d\sigma_{g_{0}}>0$. Therefore, inequality (\ref{stekloff3}) holds because of the equality and not in the strict sense.
\end{proof} 
\section*{Discussion}We showed that the scope of inequality (\ref{stekloff1}) is not as broad as the result of Proposition $2$ suggests in \cite{Escobar1}. We showed that this inequality cannot be strict since, it requires that $Q_{g}(\phi) \geq\Big(\displaystyle\max_{\partial M}e^{-f}\Big)Q_{g_{0}}(\phi)$ be strict. In fact, the latter is an equality which only occurs if the function $f$ is constant on $\partial M$; this suggests that the lower bound in Example $7$ in \cite{Escobar1} may need to be more carefully stated or perhaps revised. Similarly, more general applications of this inequality may need to be handle with care since, as it is showed, the result is not as general as it seems to be proposed. 

\section*{Declarations}
The author declares that he has no conflicts of interest.

\section*{Funding}
This work was funded by Pontificia Universidad Javeriana at Bogotá, D.C under the project ID-000000000011267.

\section*{Acknowledgment}The author would like to thank Professor Fernando Novoa Ramirez for his guidance and support.

\bibliography{REFS}
\end{document}